\documentclass[12pt]{article}
\usepackage{fullpage,graphicx,psfrag,amsmath,amsfonts,verbatim}
\usepackage[small,bf]{caption}
\usepackage{amssymb}
\graphicspath{{img/}}

\usepackage{listings}
\usepackage{longtable}
\usepackage{subcaption}

\usepackage{color}

\definecolor{codegreen}{rgb}{0,0.6,0}
\definecolor{codegray}{rgb}{0.5,0.5,0.5}
\definecolor{codepurple}{rgb}{0.58,0,0.82}
\definecolor{backcolour}{rgb}{0.95,0.95,0.98}

\lstdefinestyle{mystyle}{
backgroundcolor=\color{backcolour},
keywordstyle=\color{magenta},
numberstyle=\tiny\color{codegray},
stringstyle=\color{codepurple},
basicstyle=\ttfamily\small,
breakatwhitespace=false,
breaklines=true,
captionpos=t,
keepspaces=true,
numbers=left,
numbersep=5pt,
showspaces=false,
showstringspaces=false,
showtabs=false,
tabsize=2,
}

\definecolor{seagreen}{rgb}{0.18, 0.55, 0.34}
\definecolor{mediumviolet-red}{rgb}{0.78, 0.08, 0.52}
\definecolor{khaki}{rgb}{0.94, 0.9, 0.55}

\lstdefinelanguage{mypython}
{
keywords=[1]{from, import, as, assert, not, print, nonneg, PSD, axis},
keywordstyle=[1]{\color{mediumviolet-red}},
keywords=[2]{cp, lo, pl, cvxpy, Variable, Parameter,
	sqrt, exp, numpy, np, Problem, Minimize, Maximize, value, solve, inner,
	sum, multiply, arange, range, norm1, norm2, norm_inf, abs, square,
	diagonal, outer, pos, hstack, power},
keywordstyle=[2]{\color{seagreen}},
upquote=true,
showstringspaces=false,
basicstyle=\ttfamily,
columns=fullflexible,
keepspaces=true,
emph={True,False,def,return,float,class,match,switch,len},
emphstyle={\color{seagreen}},
belowskip=1em,
aboveskip=1em,
morecomment=[l]{\#}
}

\newcommand{\BEAS}{\begin{eqnarray*}}
\newcommand{\EEAS}{\end{eqnarray*}}
\newcommand{\BEA}{\begin{eqnarray}}
\newcommand{\EEA}{\end{eqnarray}}
\newcommand{\BEQ}{\begin{equation}}
\newcommand{\EEQ}{\end{equation}}
\newcommand{\BIT}{\begin{itemize}}
\newcommand{\EIT}{\end{itemize}}
\newcommand{\BNUM}{\begin{enumerate}}
\newcommand{\ENUM}{\end{enumerate}}

\newcommand{\BA}{\begin{array}}
\newcommand{\EA}{\end{array}}

\newcommand{\eg}{{\it e.g.}}
\newcommand{\ie}{{\it i.e.}}

\newcommand{\reals}{{\mbox{\bf R}}}

\newcommand{\symm}{{\mbox{\bf S}}}  % symmetric matrices

\newcommand{\diag}{\mathop{\bf diag}}

\newcommand{\argmin}{\mathop{\rm argmin}}

{\begin{quote}}{\end{quote}}

\makeatletter
\long\def\@makecaption#1#2{
\vskip 9pt
\begin{small}
\setbox\@tempboxa\hbox{{\bf #1:} #2}
\ifdim \wd\@tempboxa > 5.5in
\begin{center}
\begin{minipage}[t]{5.5in}
\addtolength{\baselineskip}{-0.95pt}
{\bf #1:} #2 \par
\addtolength{\baselineskip}{0.95pt}
\end{minipage}
\end{center}
\else
\hbox to\hsize{\hfil\box\@tempboxa\hfil}
\fi
\end{small}\par
}
\makeatother

\newcounter{oursection}

\newcounter{lecture}

\graphicspath{{img/}}

\usepackage{mathtools}
\usepackage{adjustbox}
\usepackage{booktabs}
\usepackage{csquotes}

\usepackage{hyperref}
\hypersetup{ colorlinks   = true,    % Colours links instead of ugly boxes
urlcolor     = blue,    % Colour for external hyperlinks
linkcolor    = blue,    % Colour of internal links
citecolor    = red      % Colour of citations
}

\title{Learning Parametric Convex Functions}
\author{Maximilian Schaller
\and Alberto Bemporad 
\and Stephen Boyd
}

\begin{document}
\maketitle

\begin{abstract}
A parametrized convex function depends on a variable and a parameter,
and is convex in the variable for any valid value of the parameter.  
Such functions can be used to specify parametrized convex 
optimization problems, \ie, a convex optimization family, in 
domain specific languages for convex optimization.
In this paper we address the problem of fitting a 
parametrized convex function that is compatible with disciplined programming,
to some given data.  
This allows us to fit a function arising in a convex optimization
formulation directly to observed or simulated data.
We demonstrate our open-source implementation on several examples,
ranging from illustrative to practical.
\end{abstract}

\clearpage
\tableofcontents
\clearpage

\section{Introduction}

\subsection{Parametrized convex functions}
A \emph{parametrized convex function} (PCF) $f$ has the form
\[
f: \reals^n \times \Theta \to \reals^d,
\]
where $\Theta \subseteq \reals^p$.
To be a PCF, $f$ must be continuous in $\theta$,
and for each $i=1, \ldots, d$, $f_i(x, \theta)$ is convex in 
$x$ for any $\theta \in \Theta$.
We refer to the first argument $x$ of the PCF $f$ as the \emph{variable},
and the second argument $\theta$ as the \emph{parameter}.
When $d=1$, we refer to $f$ as a scalar PCF.

\subsection{Disciplined convex programming}
The terms variable and parameter in a PCF 
are taken from disciplined convex programming (DCP),
a method for expressing a PCF as an expression in a domain specific language (DSL) 
constructed from variables, constants,
parameters, and a small library of functions called
atoms \cite{agrawal2018rewriting,diamond2016cvxpy}. In DCP, the 
expression must be constructed in a specific way that corresponds to a composition 
rule that establishes convexity of the function with respect to the variable, for 
any valid parameter.  

Functions expressed in DCP form can be used to form
a parametrized convex optimization problem or convex optimization family.
When the parameters are given specific numerical values, we obtain a
problem instance, which can be solved by
automatically transforming the problem instance to a canonical form, solving the canonical
form, and then retrieving the solution of the original problem instance from the solution 
of the canonicalized problem instance.
Examples of DSLs that leverage DCP for modeling convex optimization problems
(and support parameters)
include CVXPY~\cite{diamond2016cvxpy} (in Python),
CVXR~\cite{fu2017cvxr} (in R),
Convex.jl~\cite{convexjl} and JuMP~\cite{Dunning2017jump} (in Julia).
Precursors that do not handle parameters include CVX~\cite{grant2014cvx}
and YALMIP~\cite{lofberg2004yalmip} (in Matlab).

%For example
%in automatic control, also known as sequential decision making under uncertainty,
%the action in each time period is chosen 
%based on quantities that are known when the action is to be chosen
%(called the context), with
%tne mapping from the context to the action called the policy.
%A convex optimization control policy (COCP) \cite{XXX} is
%a policy in which the action is a solution of a convex optimization problem
%that is parameterized by the context.
%As a specific example, consider a trading algorithm in quantitative finance,
%which finds the action (the trade list, \ie, amounts of each asset to buy and sell),
%by solving a convex optimization problem. The variables in the problem include the
%trades and post-trade holdings; the parameters contain the context,
%which consistes of current holdings,
%forecasts of quantities such as asset returns, trade volume, and trading costs
%\cite{boyd2024markowitz,narang2013inside,markowitz1952portfolio}.

\subsection{Examples}
PCFs arise in many applications,
including control, machine learning, resource allocation, and finance, to
name just a few.  We describe below a few typical ones, some specific and some
more generic. In our examples we assume that $f$ is a scalar PCF.

\paragraph{Fuel use map.}  Here $f(x,\theta)$ gives the instantaneous fuel use rate.
The variable $x$ might correspond to thrust or power output,
variables that typically appear in an optimization problem;
the parameter $\theta$ contains additional parameters that affect the fuel use, 
such as temperature.  When modeling the total fuel use of a set of actuators
on a moving vehicle, $x$ can represent the net force and torque on the vehicle,
and $\theta$ can contain the vehicle orientation, which would affect which actuators
are used to obtain the required force and torque.

\paragraph{Battery aging model.}
Here $f(x,\theta)$ gives the rate of aging of a battery, \ie, reduction of 
capacity due to use,
relative to the initial battery capacity.
The variable $x$ corresponds to the battery charge/discharge rate,
which might appear in an optimization problem. 
The parameter $\theta$ contains other quantities
that affect aging of the battery, such as its temperature.

\paragraph{Convex optimization control policy (COCP).}
A control policy maps the context, \ie, what is known at a given time, 
such as the state of
a dynamic system and other measurable quantities, into an action denoted
$u \in \mathcal U \subseteq \reals^r$.
In a convex optimization control policy (COCP) \cite{boyd2020control} the action
$u$ is found as the solution of a 
convex optimization problem that is parametrized by the context.
As a concrete example, $f(x, \theta)$ gives the cost of the action $u$, combined
with the long term cost of the next state, which depends on the action.
The parameter $\theta$ contains the context observed when computing 
the control action
\cite{rawlings2017model, kouvaritakis2016model, garcia1989model}.

\paragraph{Resource allocation.}
Here $f(x, \theta)$ gives the cost (or negative utility) 
of providing resources specified by the vector $x \in \reals^n$ to $n$ agents.
The parameter $\theta$ might contain information like time of the year, month, or day.

\paragraph{Financial portfolio construction.}
The goal of financial portfolio construction is to find a portfolio of
financial assets that maximizes the expected return of the portfolio while limiting
the investment risk. The function $f(x, \theta)$ might represent a combination of 
expected return and risk as a function of $x$, which describes the portfolio.
The parameter $\theta$ provides information about market conditions, or current 
forecasts.

\subsection{Learning a parametrized convex function}
This paper concerns learning a PCF $f$ from data, \ie, 
fitting a PCF to data,
\BEQ\label{e-loss+reg}
(x^k, \theta^k) \in \reals^n \times \Theta, \quad y^k \in \reals^d, 
\quad k = 1,\ldots,N.
\EEQ
The PCF $f$ is specified by its architecture and a choice of \emph{model weights},
which we denote $w \in \reals^q$.
We require our approximation $f$ to be DCP expressible, which implies that it can be used
to construct parametrized convex problems.
This fitting method allows us to learn a DCP expression directly from observed 
(or generated) data.

\paragraph{Fitting method.}
We use a standard regularized loss method to choose a set of model weights.
Let $\ell: \reals^q \times \reals^n \times \Theta \times \reals^d \to 
\reals$ denote a loss function, 
and $r: \reals^q \to \reals$ a regularizer function.
The loss function is used to judge how well our approximation fits the data,
and the regularization function is meant to penalize complexity.
We choose $w$ to (approximately) minimize the regularized average loss,
\BEQ
\frac{1}{N} \sum_{k=1}^N \ell(w; x^k, \theta^k, y^k) + \lambda r(w),
\label{eq:training_loss}
\EEQ
where $\lambda \geq 0$ is a hyper-parameter that scales the regularization.
In this context we refer to the data as training data, since it is used to train
or learn the PCF.

As in all data fitting methods we judge a candidate PCF by its accuracy on unseen data.
(The loss used to evaluate the model 
can differ from the loss $\ell$ used to train the model.)
To find a PCF that performs well on unseen data we use the standard technique
of cross-validation.  We partition the original data into $K$ groups of 
approximately equal size, called folds, and for each fold we train the
model weights on the data not including that fold, and evaluate the average loss 
on the data in that fold. 
We average the losses for the $K$ folds to obtain an overall fitting metric.
We choose the architecture and the regularization hyper-parameter $\lambda$
so as to minimize the overall fitting metric. Once the architecture and 
hyper-parameter are chosen, we fit the PCF using all the data.

\paragraph{Special cases.}
The PCF learning problem, \ie, minimizing \eqref{e-loss+reg},
reduces to well known problems in special cases.
When $n=0$, \ie, there is no variable, it reduces to the very general 
problem of fitting a continuous function from $\Theta$ into $\reals^d$.
When $p=0$ (or $\Theta$ is a singleton), there is (effectively) no parameter, 
and the problem reduces
to fitting a convex function to some given data.  Several architectures have
been developed for this task, such as input-convex neural networks
\cite{amos2017input, deschatre2025input, liu2025icnn}, which we build on.

\paragraph{Non-parametric PCF fitting problems.}
We can solve several more general fitting problems exactly, with no constraint
on the architecture (\ie, non-parametric) and no regularization.
When $p=0$ and the loss is convex in $y$, the problem of minimizing the average loss
over all convex functions, which is an infinite-dimensional non-parametric fitting 
problem, can be solved exactly using convex optimization,
as decribed in \cite[\S5.5.5]{boyd2004convex}.

This can be extended to the more general case when $p>0$ and $r=0$.
To do this we collect the data into groups associated with unique values of 
$\theta$, and 
for each one, we fit a convex function to these data as described above.  We can 
then interpolate these functions for values of $\theta$ not appearing in the data.
This solution globally minimizes the loss, but since it does not 
include any penalty on complexity of $f$,
it is very likely to perform poorly on unseen data.

\subsection{Contribution}

We provide a seamless path from data to a PCF that
can be used for parametrized convex optimization.
Our open-source implementation LPCF offers a simple user 
interface for fitting a PCF to variable-parameter-output triples
and a variety of extensions for PCFs with special properties and uses.
The resulting PCFs can be exported for immediate use in optimization
frameworks like JAX~\cite{jax2018github} or CVXPY.

\subsection{Outline}

In \S\ref{s-architecture} we propose an architecture,
which is a generalization of input-convex neural networks that handles
parameters, and describe our specific implementation.
In \S\ref{s-extensions} we describe a set of extensions that allow for modeling
more specialized types of PCFs.
In \S\ref{s-experiment} we give some numerical results for both
simple illustrative examples and some practical ones.
In \S\ref{s-conclusions} we conclude the paper.

\section{Proposed neural network architecture}\label{s-architecture}

\subsection{Architecture}
Let $\phi:\reals \to \reals$ be a nondecreasing convex activation function,
such as rectified linear unit (ReLU) with 
$\phi(a)= \max\{a,0\}$ or softplus with $\phi(a)=\log(1+e^a)$,
that we will also refer to as \emph{logistic} due
to its use in formulating logistic regression problems.
We have $L$ layers, with $z^l \in \reals^{n_l}$ the activation of 
layer $l$, $l=1,\ldots, L-1$.
We have
\BEQ\label{e-arch}
z^0 = x, \quad
z^l = \phi \left(W^l z^{l-1} + V^l x + \omega^l\right), \quad l=1,\ldots,L-1, \quad
y=W^L z^{L-1} + V^Lx + \omega^L,
\EEQ
where $\phi$ is applied componentwise.
Here $W^l \in \reals^{n_l \times n_{l-1}}$
is the weight matrix associated with layer $l$, 
$V^l \in \reals^{n_l \times n}$ is the
weight matrix that feeds the input $x$ into layer $l$, and 
$\omega_l \in \reals^{n_l}$ is the offset for layer $l$.
%We have $n_1=n$ and $n_L=d$.
We can take $W^1 = 0$ without loss of generality since $z^0=x$.
This is a standard residual network architecture, with feedforward from
the input $x$ into each layer \cite{he2016deep,he2016identity}.

We take the weight matrices $W^l$, $V^l$, and offsets $\omega^l$ to be a
function of the parameter $\theta$,
\BEQ\label{e-param-network}
(W^2, \ldots, W^L, V^1, \ldots, V^L, \omega^1, \ldots, \omega^L) = \psi (\theta),
\EEQ
where $\psi:\Theta \to \reals^m$ describes a generic neural network
architecture with output dimension
\[
m = n_2 n_1 + \cdots + n_L n_{L-1} + n_1 n + \cdots + n_L n + n_1 + \cdots + n_L.
\]
(We will impose one constraint on the function $\psi$, described below.)
We denote the weight matrices and offsets as $W^l(\theta)$, $V^l(\theta)$ and
$\omega^l(\theta)$ to emphasize their dependence on the parameter $\theta$.
The model weights $w$ defining the loss in~\eqref{eq:training_loss}
are the weights defining $\psi$.

\begin{figure}
\centering
\includegraphics[width=0.9\columnwidth]{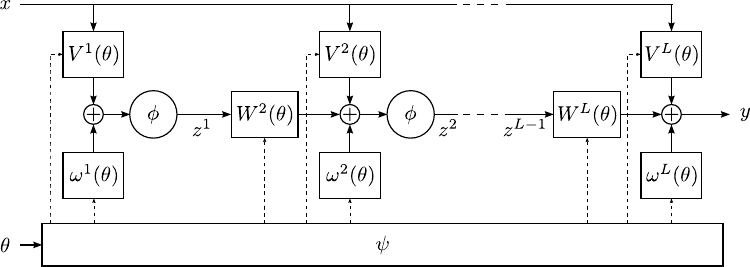}
\caption{Neural network architecture for PCF $y=f(x, \theta)$.}
\label{fig:architecture}
\end{figure}

The proposed architecture is given by \eqref{e-arch} 
and \eqref{e-param-network}.
It is specified by the layer dimensions $n_1, \ldots, n_{L-1}$
(the final layer width $n_L$ is fixed to $d$) and 
the architecture for $\psi$.
This architecture determines a function $f: \reals^n \times \Theta \to \reals^d$,
illustrated as a block diagram in figure \ref{fig:architecture}.
The model weights $w$ associated with $f$ are the model weights appearing in 
the righthand side of \eqref{e-param-network}.

\paragraph{Convexity.}
We impose one restriction on the weight matrices: $W^l(\theta)$ are 
elementwise nonnegative, for any $\theta \in \Theta$.
This can be enforced by the architecture of $\psi$, for 
example by having the weight matrices $W^l(\theta)$ come directly from a
nonnegative activation function such as ReLU or logistic,
without any offset.

With this restriction, the function $f$ is a PCF.
We argue using recursion that each element of $z^l$ is a convex
function of $x$, for any fixed $\theta$.
It is evidently true for $l=1$. 
Assuming that each element of $z^{l-1}$ is a convex function of $x$,
we observe that each entry of
\[
W^l(\theta) z^{l-1} + V^l(\theta) x + \omega^l(\theta)
\]
is a convex function of $x$, since it is a nonnegative weighted sum of 
convex functions, plus an affine function of $x$.
By the composition rule \cite[\S 3.2.4]{boyd2004convex}, 
each entry of $z^l$ is a convex function of $x$.
This argument is exactly the one used in DCP. This means that, assuming the 
activation function is an atom, $f(x,\theta)$ is DCP.

As in an input convex neural network, our architecture (and the nonnegativity
constraint on $W^l$), was chosen 
specifically so that it implements a PCF for any valid choice of
the weights. 
%We note that the parameter $\theta$ enters the network that maps $x$ to $y$ only 
%via the offsets at each layer.  A more expressive architecture allows 
%the weight matrices $V^l$ and $W^l$ to also depend on $\theta$, 
%in addition to the offsets (assuming the nonnegativity constraint on $W^l$
%is respected).  Thus also results in a PCF, but it is not DPP, since
%the way $V^l$ and $W^l$ enter $f$ violates the DPP rules.

\paragraph{Loss and regularizer.}
The loss and regularizer used in the fitting process are arbitrary.
Some conventional choices of loss include quadratic, $\ell_1$, or
Huber \cite{huber1992robust, huber2011robust}.
Conventional choices of regularizer include quadratic, $\ell_1$, or a combination
(known as elastic net regularization).
See, \eg, \cite{zou2005regularization} or \cite[Chap.~6.3]{boyd2004convex}.

\subsection{Implementation}\label{s-implementation}

Our open-source implementation \verb|lpcf| is available at
\begin{center}
\url{https://github.com/cvxgrp/lpcf}.
\end{center}
In \verb|lpcf|
the user can specify the architecture, activation function, and type of 
regularization. The generic \verb|.fit(data)| method fits the PCF to data, using
cross-validation to choose the regularization hyper-parameter, and reports
the performance of the fit obtained, for separate test data.
In \verb|lpcf|, the R2-score~\cite{draper1998applied,pearson1905general}
is the standard validation metric used for cross-validation and testing.
Utilities are provided for the user to evaluate the resulting PCF, score it 
on another data set, and to export the function to CVXPY, where it can be 
freely used anywhere a convex function can be.

A simple script illustrating this functionality is shown in 
figure~\ref{code:export_cvxpy}.
In lines 5--6 we instantiate a \verb|pcf| object and fit its weights to data.
In line 11 we export to the DCP expression \verb|f| via the \verb|tocvxpy|
method, which takes a CVXPY variable and a CVXPY parameter as arguments,
which are defined in the two lines above.
We illustrate the use of \verb|f| in constructing a CVXPY problem
in line 15.
In this example we simply add the fitted PCF to the objective, but we note
that it can be used anywhere in CVXPY that a convex function appear, \eg,
in constraints.

\begin{figure}
	\lstset{language=Python,
		numbers=left,
		xleftmargin=0.07\columnwidth,
		linewidth=\columnwidth}
	\begin{lstlisting}[frame=lines]
from lpcf.pcf import PCF
import cvxpy as cp

# create default PCF object and fit to data Y, X, Theta
pcf = PCF()
pcf.fit(Y, X, Theta)

# export to cvxpy expression
x = cp.Variable((n, 1))
theta = cp.Parameter((p, 1))
f = pcf.tocvxpy(x, theta)

# solve cvxpy problem involving f
# g is another function, constraints a list of (in)equalities
problem = cp.Problem(cp.Minimize(f + g), constraints)
theta.value = ...
problem.solve()
	\end{lstlisting}
	\caption{Using LPCF with CVXPY. The initialization code for dimensions
	\texttt{n}, \texttt{p}, data \texttt{Y}, \texttt{X}, \texttt{Theta},
	and CVXPY objects \texttt{g} and \texttt{constraints} is omitted for clarity.}
	\label{code:export_cvxpy}
\end{figure}

To fit the weights $w$ of the network $\psi$ we use \verb|jax-sysid|~\cite{Bem24},
a Python package based on the auto-differentiation framework JAX,
for system identification, neural network training, and nonlinear regression/classification.
The package supports several types of regularization, simple bounds on $w$,
parallel training from different initial values of $w$, and quasi-Newton methods for
faster terminal convergence and better model quality than gradient descent.

\begin{figure}
\centering
\includegraphics[width=0.9\columnwidth]{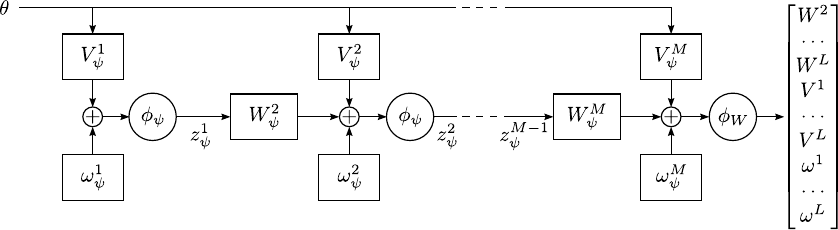}
\caption{Neural network architecture for $\psi(\theta)$.
%\AB{The output activation $\psi_+(\cdot)$ that makes $W^l\geq 0$ is missing,
%maybe we should add it.}
}
\label{fig:psi}
\end{figure}

\paragraph{Default choices.}
We make a number of default choices, all of which can be modified by the user.
We take $L=3$ layers in the input-convex network, with layer dimensions
$n_1 = \cdots = n_{L-1} = 2 \lfloor (n + d) / 2 \rfloor$,
and ReLU activation function.
We take $\psi$ to be a fully connected neural network with feedforward terms
from $\theta$ into each layer, very similar to the architecture of the input-convex
network. The output activation $\phi_W$ only affects $W^2, \ldots, W^L$.
We use the ReLU function to make them elementwise nonnegative.
The architecture is visualized in figure~\ref{fig:psi}.
By default there are $M=3$ layers, 
of which the two inner ones are $\lfloor (p + m) / 2 \rfloor$ wide, with 
ReLU activation.
By default, we use a quadratic loss (mean squared error) without regularization,
\ie, $\lambda =0$, and the cross-validation
procedure is turned off.  Otherwise, the default number of folds for
cross-validation is 5.
The default optimizer runs 200 iterations of Adam~\cite{kingma2014adam} to obtain
a good set of network weights,
followed by 2000 iterations of
L-BFGS-B~\cite{byrd1995limited} 
(a variant of L-BFGS for bound-constrained optimization,
to respect the nonnegativity constraints on $W^l$) to refine the model weights.
We run the entire fitting method on multiple random initial sets of model weights,
and take the best one as our final choice.
We parallelize training with multi-core processing, with a default of 4 cores.
By default, we set the number of initializations to 10 or
the number of cores, whichever is greater.
All of these settings can be customized by the user.

\section{Extensions}
\label{s-extensions}

We describe a few extensions of the basic model and methods.

\subsection{Adding a quadratic term}
\label{s-quadratic}
We can add a (convex) quadratic function to our more general
architecture~\eqref{e-arch}, 
\[
y=x^T Q x + W^L z^{L-1} + V^Lx + \omega^L,
\]
where $Q \in \symm_+^n$.  
We use the representation $Q = U^T U$ with upper-triangular
$U \in \reals^{n \times n}$, and emit $U$ from $\psi(\theta)$, as in
\[
(W^2, \ldots, W^L, V^1, \ldots, V^L, \omega^1, \ldots, \omega^L,  U) = \psi (\theta).
\]
This automatically guarantees positive semidefiniteness of $Q$, without any additional
constraints on the parameter network $\psi$.

A further variation is a low rank plus diagonal quadratic
$Q = F^T F + \diag (d_1^2, \ldots, d_n^2)$,
where $F \in \reals^{m \times n}$ is a wide matrix with $m \ll n$ and $d \in \reals^n$
is a vector of diagonal entries before they are squared to make 
$Q$ positive semidefinite.
Such a form can be used to capture the dominant directions of curvature when $n$ is large.

\subsection{Monotonicity}
\label{s-monotonicity}
In certain cases we may wish to impose that the PCF $f$ is monotonically
nondecreasing or nonincreasing with respect to some or all the components of $x$
for any $\theta\in\Theta$, \eg, when modeling a concave and 
increasing utility function (that we wish to maximize).

Monotonicity of $f$ can be imposed as follows. Assume that the activation function
$\phi$ is nondecreasing (\eg, ReLU or softplus).
Consider the architecture of the network $\psi$ described in figure~\ref{fig:psi}.
We achieve monotonicity by obtaining $V^j$, $j=1,\ldots,L$ similary to how we obtain $W^j$,
as the output of the nonnegative activation $\phi_W$, $j=2,\ldots,L$
(enforcing convexity in that case).

We can prove by induction that the resulting function $f$ is monotonically increasing
with respect to each component $x_i$, $i=1,\ldots,n$, for any $\theta\in\Theta$.
First, $z^1=\phi(W^1x+V^1x+\omega^1)=\phi(V^1x+\omega^1)$ is increasing,
since $\phi$ is increasing and its argument is affine and increasing with respect to each $x_i$;
assuming that $z^l=\phi(W^lz^{l-1}+V^lx+\omega^l)$ is increasing in $x_i$,
then $W^{l+1}z^l+V^{l+1}x+\omega^{l+1}$ remains increasing
since both $W^{l+1}$ and $V^{l+1}$ are nonnegative.
Since $\phi$ is increasing, it follows that $z^{l+1}$ is increasing in $x_i$.

By cascading instead $V^j$ with a nonpositive activation function $\phi_-$,
such as $\phi_- = -\phi_+$, $j=1,\ldots,L$, by repeating the above argument it is easy to see
that the resulting PCF $f$ is convex and decreasing with respect to
each component of $x$ for any $\theta\in\Theta$.

\subsection{Specifying a subgradient} \label{s-argmin}
We may wish to specify (or encourage) a PCF to be minimized
at a particular point $g(\theta)$,
where $g:\Theta \to\reals^n$ is given.
When $f$ is differentiable this is equivalent to $\nabla_x f(g(\theta),\theta)=0$, or,
more generally, that $0\in \partial_x f(g(\theta),\theta)$ (the subdifferential
of $f$ with respect to $x$).
This can be encouraged in the training process by adding a regularization term 
such as
\BEQ
\ell_{\rm min}(w) = \frac{\rho_{\rm min}}{N}\sum_{i=1}^N \left\|\nabla_x f(g(\theta^k),\theta^k)\right\|_2^2
\label{eq:argmin_loss}
\EEQ
to the loss~\eqref{eq:training_loss}, where
$\rho_{\rm min}$ is a hyper-parameter
that scales the regularization. This (sub)gradient can be computed by automatic 
differentiation of the PCF $f$ with respect to its first argument.

As a further generalization, we can require (or encourage) $f(x, \theta)-q(\theta)^Tx$ 
to be minimized at $x=g(\theta)$.  This is equivalent to 
$\nabla_x f(g(\theta),\theta)=q(\theta)$, and a similar 
regularization term can be added to the training objective.

\subsection{Fitting a parametrized convex set}
\label{s-classification}

The same framework can be used to fit a parametrized convex set $C:\Theta \to
2^{\textbf{R}^n}$, described as
\[
C(\theta) = \{ x \mid f(x,\theta) \leq 0 \},
\]
where $f$ is a PCF. 
Our data has the form
$(x^k, \theta^k, y^k)$, $k=1, \ldots, N$,
with $y^k \in \{-1, 1\}$,
where $y^k = -1 $ means that $x^k \not\in C(\theta^k)$ and $y^k=1$ means that
$x^k\in C(\theta^k)$.
%As a simple extension of standard linear classification, we call it
%\emph{convex classification} when the PCF $f$ is used to model the probability
%\[
    %\prob(y^k=1) = \frac{1}{1+e^{-f(x^k,\theta^k)}},
%\]
%so that $f(x^k,\theta^k)\leq0$ implies $\prob(y^k=-1) \geq 1/2$ and $f(x^k,\theta^k) \geq 0$
%implies $\prob(y^k=1) \geq 1/2$. The logarithm of the likelihood function 
%\[
%\begin{aligned}
%{\mathcal L}(w|Y)&=\prod_{k:\ y_k=1}\prob(y^k=1)\prod_{k:\ y_k=-1}\prob(y^k=-1)\\
%&=\prod_{k:\ y_k=1}\frac{1}{1+e^{-f(x^k,\theta^k)}}
%\prod_{k:\ y_k=-1}\left(1-\frac{1}{1+e^{-f(x^k,\theta^k)}}\right)=\prod_{k=1}^N\frac{1}{1+e^{-y^k f(x^k,\theta^k)}}
%\end{aligned}
%\]
%associated with the sequence $Y=(y_1,\ldots,y_N)$ is maximized by
%defining in~\eqref{eq:training_loss} the loss
We use the logistic loss function
\[
\ell(w; x^k, \theta^k, y^k)=\log (1+e^{-y^k f(x^k,\theta^k)}).
\] 
We can judge the loss on test data using either the same logistic loss function
or the actual error rate, 
\[
\ell^\text{test}(w; x^k, \theta^k, y^k)= \left\{ \begin{array}{ll} 
0 & y^k f(x^k,\theta^k) \geq 0\\
1 & y^k f(x^k,\theta^k) < 0.
\end{array}\right.
\] 

\section{Experiments}\label{s-experiment}

We start with two illustrative examples on simple functions,
and then give two real applications, battery aging and control.
We conduct the experiments on an Apple M4 Max machine.

\subsection{Piecewise affine function on $\reals$}

We generate data from a piecewise affine (PWA) function of the form
\[
f^\text{true}(x, \theta) = s_+ \max\{0,x - m\} + s_-\max\{0,m - x\} + v,
\]
where $x\in \reals$ and $\theta = ( s_+, s_- , m, v ) \in \reals^4$.
Note that $f^\text{true}$ is convex (\ie, a PCF) when $s_+ \geq - s_-$, 
but we also carry out experiments when this is not the case.
We use a quadratic loss function.

\paragraph{Experimental setup.}
We sample 2000 values of $\theta$ from a uniform distribution on $[-1,1]^4$.
For each value of $\theta$, we take 50 equally spaced
points $x$ on $[-1,1]$, which gives $N = 2000 \times 50 = 10^5$ data points.
We divide these data into two sets: Those for which $f^\text{true}$ is 
convex, and those for which it is not.
For those that are nonconvex, we find the best affine fit, which is
the (non-parametric) closest convex approximation of the function.

We fit this data using the default values for the network architecture and learning
parameters.

\paragraph{Results.}
Table~\ref{tab:pwa} contains the root mean squared error (RMSE) values of the
learned PCF $f$, for $10^5$ random test data points with the true function values
on a scale from $-3$ to $3$.
The first three values are with respect to the true
data-generating function $f^\text{true}$, for all test data and computed for
the convex and nonconvex data separately. The last value is the RMSE for
nonconvex data, computed with the best affine fit as the ground truth.
\begin{table}
\centering
\begin{tabular}{l|r}
Data & RMSE \\ \hline 
All & 0.054\\
Convex & 0.001\\
Nonconvex & 0.077\\
Nonconvex (relative to best affine fit) & 0.004\\
\hline
\end{tabular}
\caption{RMSE values. The value of $f^\text{true}$ ranges between $-3$ and $3$.
The last entry compares the PCF approximation for nonconvex data to the
best affine fit.}
\label{tab:pwa}
\end{table}
We observe that the overall RMSE is very small. When only considering the convex examples,
the RMSE is more than an order of magnitude smaller. As expected, the RMSE is considerably
larger when only taking the nonconvex examples.
When computed with respect to the best affine fit $f^\text{linear}$,
the RMSE for the nonconvex examples is also in the order of the RMSE for the convex data.

\begin{figure}
\centering
\includegraphics[width=\columnwidth]{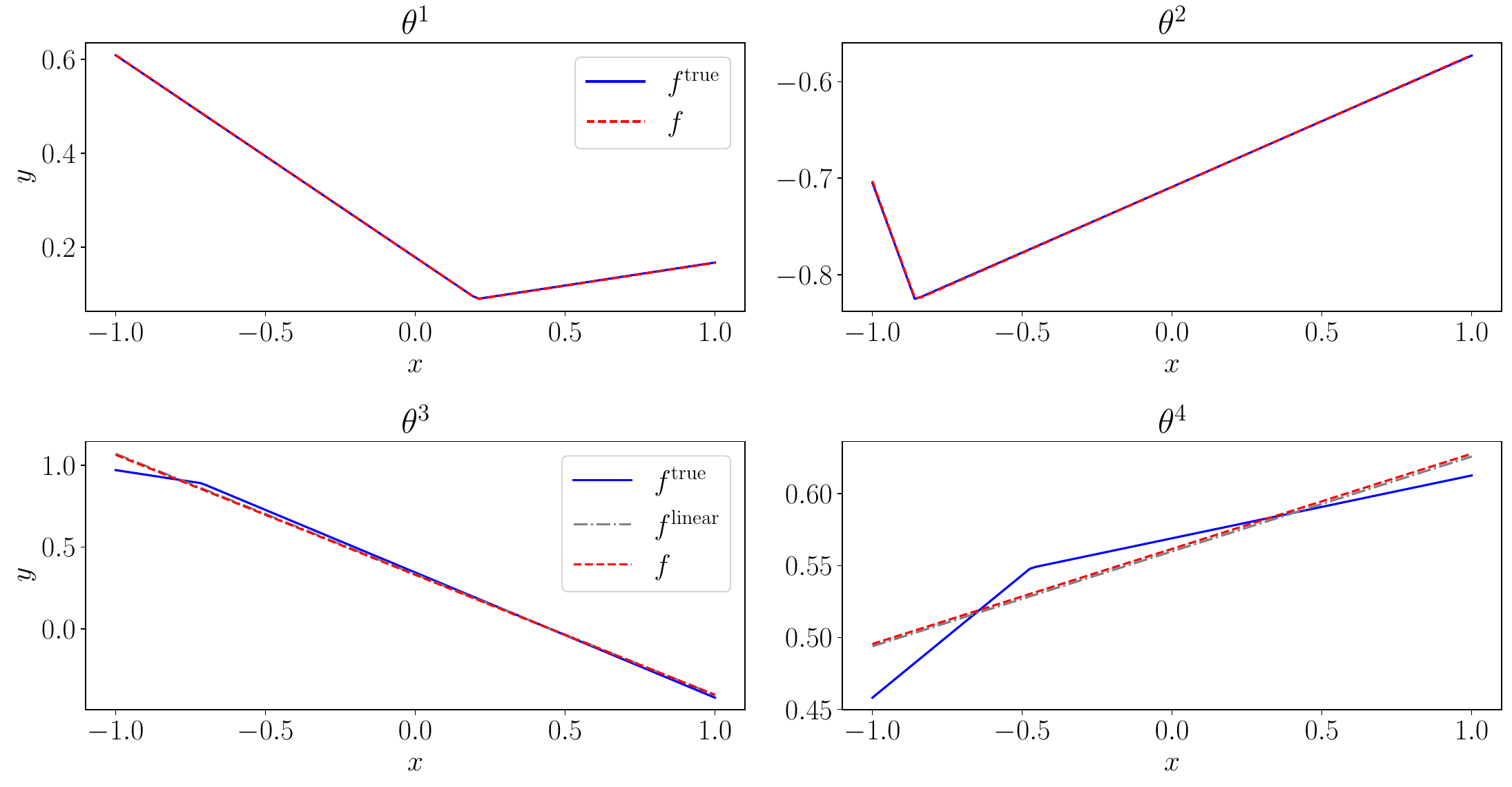}
\caption{Data-generating function $f^\text{true}$ and learned PCF $f$ for 
four parameter values. \emph{Top.} Convex $f^\text{true}$. \emph{Bottom.} 
Non-convex $f^\text{true}$.
}
\label{fig:pwa}
\end{figure}

Figure~\ref{fig:pwa} shows the true data-generating function $f^\text{true}$ and the
learned PCF $f$, for four random values of $\theta$.
For the first two parameter values $\theta^1$ and $\theta^2$, 
$f^\text{true}$ is convex; for the last two parameter values $\theta^3$ and $\theta^4$,
it is not convex.
For the nonconvex pair we also show the best affine fit.
We can see that the approximations are very good in both cases.

\clearpage

\subsection{Quadratic function on $\reals^3$}

We generate data from a parametrized quadratic function
\[
f^\text{true}(x, \theta) = x^T \theta x,
\]
where the variable is $x \in \reals^n$ and the parameter is $\theta \in \symm^n_+$
(the set of $n \times n$ positive semidefinite symmetric matrices).
The true function $f^\text{true}$ is convex, \ie, a PCF.

\paragraph{Experimental setup.}
We take $n=3$ and generate $1000$ values of $\theta$ where each entry is
taken from $[-1, 1]$.
For each value of $\theta$, we sample $100$ values of $x$ from a uniform distribution
on the unit ball, which gives $N=10^5$ data points.

We fit the PCF with the softplus activation function
and otherwise default settings for the network architecture and learning algorithm.

\paragraph{Results.}
For $10^5$ random test data points, where the value of the true data generating function $f^\text{true}$
moves between $0$ and $2$, we attain a low RMSE value of about $0.02$.

\clearpage

\subsection{Battery aging}
When repeatedly charging and discharging the battery of, \eg, a hybrid electric vehicle
or an industrial energy storage system,
the battery capacity degrades over time, a process called
battery aging~\cite{suri2016control, serrao2011optimal, ebbesen2012battery}.
We denote by $y$ the aging rate and
generate data from the battery aging model as introduced in~\cite{suri2016control}
and used in~\cite{nnorom2025aging} for optimal battery management,
\[
f^\text{true}(x, \theta) = z A^{z-1} b (\alpha q/Q + \beta)
\exp \left(\frac{-E_a + \eta b / Q}{R_g (T_0 + T)}\right),
\]
where $x = (q, b)\in \reals_+^2$ is the variable, consisting of the battery's
charge $q$ and the absolute charge rate $b$.
The parameter is $\theta = (A, Q, T) \in \reals_+^3$ and contains the
accumulated charge throughput $A$, the battery capacity $Q$, and temperature $T$.
The remaining symbols are constants. We use the values in \cite{nnorom2025aging},
the physical constants
\[
E_a = 31500,  \quad R_g=8.3145, \quad T_0=273.15,
\]
and battery parameters
\[
\alpha=28.966, \quad \beta=74.112, \quad z=0.6, \quad
\eta=152.5.
\]

\paragraph{Experimental setup.}
We generate $1000$ values of $\theta$, 
with accumulated charge throughput $A \in [0,50]$ and temperature $T \in [10, 50]$.
We keep the battery capacity fixed at $Q = 1$, since it is the slowest changing parameter
(especially for a new battery where $A$ moves fast) and one can re-fit the PCF as
$Q$ changes over time.
For each value of $\theta$, we sample $100$ values of $x$
with state of charge $q \in [0.2, 0.8]$ and charge rate $b \in [0, 30]$,
which gives $N = 10^5$ data points.

We fit the PCF with the softplus activation function.  The input-convex 
network is 5 wide,
the network $\psi$ has a single hidden layer of width 10, and we train with 
1000 and 4000 epochs with Adam and L-BFGS-B, respectively.

We compare our fit to the convex approximation used in~\cite{nnorom2025aging}
for short term battery management,
\[
f^\text{short} (x, \theta) = \mu (1 + \nu Q/2) b,
\]
where
\[
\mu = \beta \exp \left( \frac{-E_a}{R_g (T_0 + T)} \right) z A^{z-1}, \quad
\nu = \frac{\alpha}{\beta Q}.
\]
This is a PCF itself, and can be derived as the first-order Taylor expansion
of $f^\text{true}$ around the point $(q, b) = (Q/2, 0)$~\cite{nnorom2025aging}.

\begin{figure}
\centering
\includegraphics[width=\columnwidth]{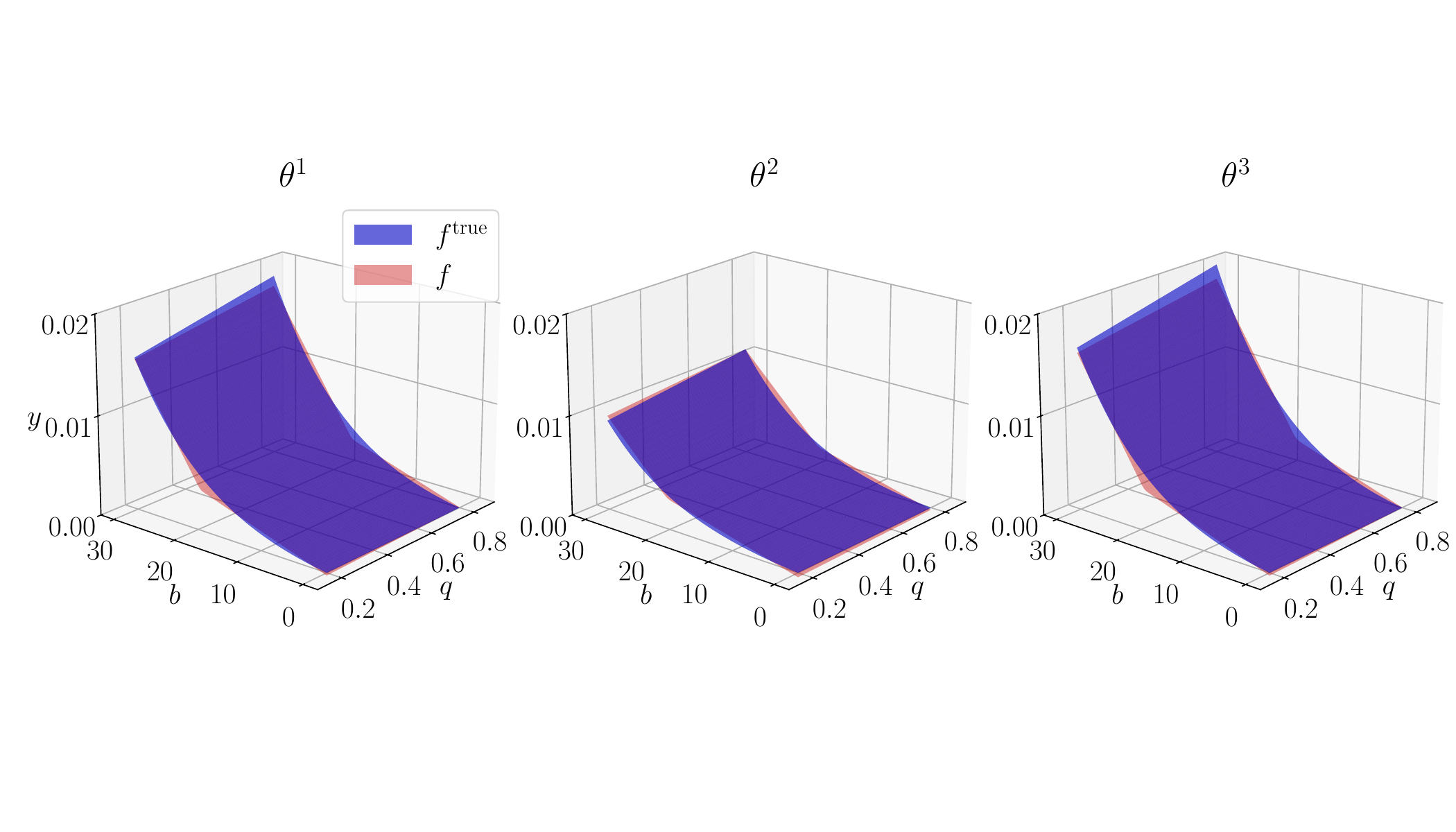}
\caption{Data-generating function $f^\text{true}$ and learned PCF $f$ for 
three parameter values.}
\label{fig:battery}
\end{figure}

\paragraph{Results.}
For $10^5$ random test points with $f^\text{true}$ on a scale from $0$ to $0.06$,
the RMSE is $0.001$.
This is about 7 times better than using the short term approximation $f^\text{short}$,
which gives an RMSE of $0.007$.
Figure~\ref{fig:battery} shows the true data-generating function $f^\text{true}$ and the
learned PCF $f$, for three random values of $\theta$.
We observe that the approximations are good in all cases.

\clearpage

\subsection{Approximate dynamic programming}
We consider a nonlinear dynamical system of the form
\BEQ
z_{t+1}= F(z_t,\theta)+G(z_t,\theta)u_t, \quad t=0,1,\ldots ,
\label{eq:affine-dynamics}
\EEQ
where $z_t\in\reals^{n}$ is the state, $u_t\in\reals^{m}$ is the input, 
$\theta$ is a vector of parameters,
$F:\reals^{n} \times \Theta \to\reals^{n}$ gives the dynamics,
and $G:\reals^{n} \times \Theta \to\reals^{n\times m}$ gives the input-to-state
matrix.
Since $z_{t+1}$ is an affine function of $u_t$, this is called input-affine form.
We are given the inital state $z_0$, and seek inputs $u_0, u_1, \ldots $ that minimize
the cost function
\BEQ
J(z_0) = \sum_{t=0}^{\infty}H(z_{t},u_t,\theta),
\label{eq:optimal_control}
\EEQ
where $H: \reals^n \times \reals^ m \times \Theta \to \reals_+$ 
is a nonnegative stage cost, assumed convex with respect to $(z_{t},u_t)$.
This is a nonconvex optimal control problem, which is hard to solve globally.
A local method can however be used to approximately solve this problem, typically out to
some large terminal time $t=T$.

Approximate dynamic programming (ADP) is a heuristic for approximately solving
the optimal control problem, choosing $u_t$ one step at a time by solving a 
convex problem. It has the form
\BEQ\label{e-adp}
u_t = \argmin_u \left( H(z_t,u,\theta) + 
\hat V(F(z_t,\theta)+G(z_t,\theta) u,\theta) \right),
\EEQ
for $t=0,1,\ldots $, where we update $z_t$ using the dynamics 
\eqref{eq:affine-dynamics}.
Here $\hat V: \reals^n\times \Theta \to \reals$ is a PCF, so the minimization
over $u$ is a convex optimization problem.

ADP is motivated by the Bellman or recursive form of the optimal input sequence $u_t$,
given by \eqref{e-adp} with $\hat V$ replaced by the value or cost-to-go function
\[
V(z_0,\theta) = \min_{u_0, u_1, \ldots} J(z_0,\theta)
\]
(see, \eg, \cite{bertsekas2005dynamic, stellato2017thesis,AZB24}). 

To find the convex approximate value function, we use a local method 
to approximately solve the problem (out to some large time period $t=T$) 
for multiple values of the initial state $z_0$.  Each such optimization in fact
gives us a set of (approximate) values of the value function, one obtained
at each state on the trajectory, with value equal to the corresponding tail cost.

\paragraph{Experimental setup.} 
We consider the problem of swinging up a pendulum to the vertical position 
by controlling the applied torque. We start with a nonlinear model of the
continuous time dynamics
\[
    ml^2\ddot\delta+b\dot\delta+mgl\sin \delta=u,
\]
where $\delta$ is the angular position of the pendulum, $m \in [0.5,2]$ is its mass,
$l=1$ its length, $b=0.05$ is the damping coefficient, $g=9.81$ is the gravitational
acceleration, and $u$ is the applied torque, all in standard metric units.
We obtain an input-affine discrete-time model as in~\eqref{eq:affine-dynamics} 
by setting $z=(\delta, \dot\delta)$ and using a first-order forward Euler integration 
(with sampling time $0.02$). 
The stage cost is
\[
H (z,u) = (\delta-\pi)^2+ 0.01 \dot\delta^2+ 0.001u^2.
\]
Our parameter is $\theta = m$, with $\Theta = [0.5,2]$.

We generate $N=1000$ data points by sampling uniformly 
$\delta_k\in[-\pi / 6, 7\pi / 6]$,
$\dot\delta_k\in[-1,1]$, and $m_k\in[0.5,2]$.
For each combination, we (approximately) solve the optimal 
control problem~\eqref{eq:optimal_control} out to $t=T=150$,
using the L-BFGS-B optimizer. 

We fit the PCF with the softplus activation function,
width 20 for the input-convex network,
two hidden layers of width 10 for the network $\psi$,
and default architecture otherwise.
%We select a PCF architecture with $L=3$ layers with $n_1=n_2=20$ hidden units plus a quadratic function, a neural network $\psi$ with two hidden layers with $10$ units each and $546$ outputs;
%both networks use the softplus activation function. The total number of model weights $w$ is $q=6692$.
%\MS{Do we need to give the output dimension of $\psi$ and the total number of model weights?
%I think all other quantities in this paragraph are chosen quantities and these two follow.}
We use the elastic net regularization term $r(w)=10^{-8}\|w\|_2^2+0.1\|w\|_1$
in~\eqref{eq:training_loss}
and two of the extensions described in~\S\ref{s-extensions}: We add a quadratic
term (see~\S\ref{s-quadratic}) and specify a subgradient (see~\S\ref{s-argmin}).
In particular, we know that $V(z,\theta)$ has a minimum at the equilibrium
$z^{\rm eq}=(\pi, 0)$. We promote that the same holds for $f$,
by using $g(\theta) = z^{\rm eq}$ in the regularization
term~\eqref{eq:argmin_loss}.

We train the PCF from 16 different initial values, running $1000$ Adam iterations
followed by up to $5000$ function evaluations in L-BFGS-B.

\begin{figure}
\centering
\includegraphics[width=0.45\columnwidth]{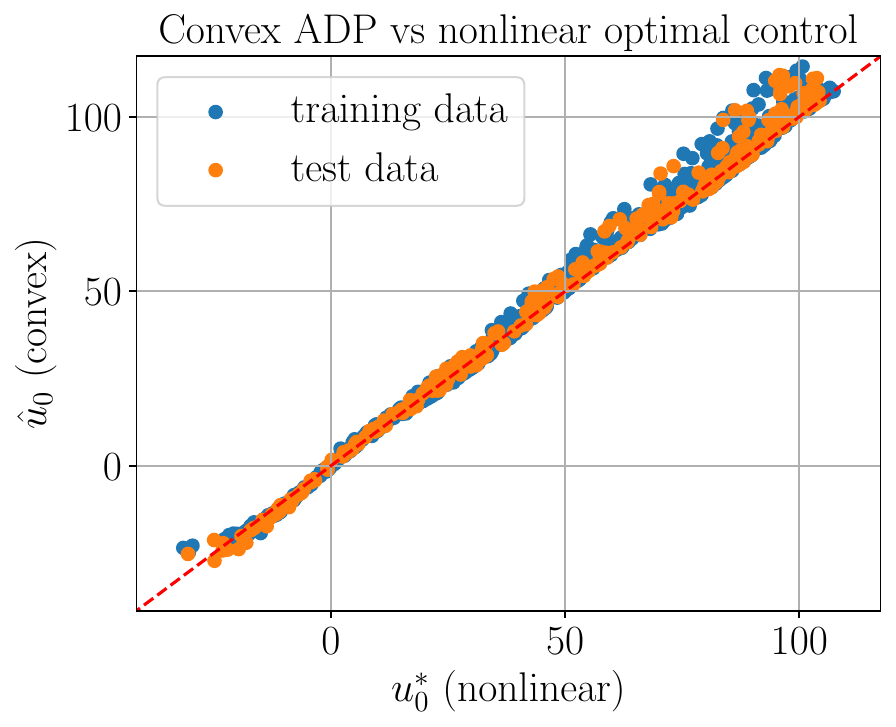}
\includegraphics[width=0.45\columnwidth]{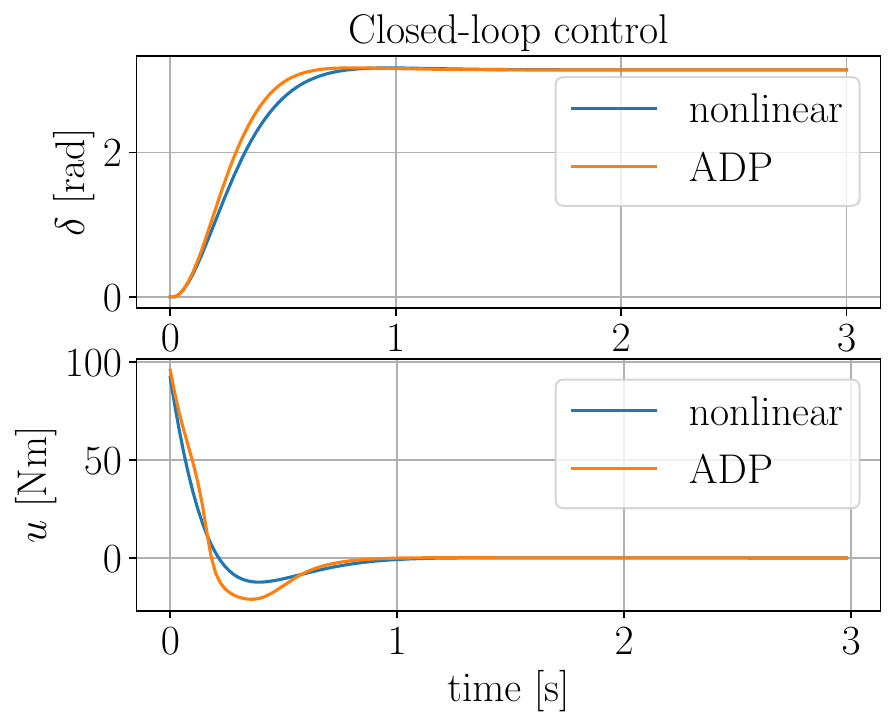}
\caption{Left: ADP control input $\hat u_0$ versus nonlinear control input $u^\star_0$.
Right: Simulation of the two controllers starting from $z_0=(0, 0)$, where $\theta = m=1$.}
\label{fig:adp}
\end{figure}

\paragraph{Results.}
%The total training time is about $3.5$ minutes.
%We select the model with the highest $R^2$-score on training data, which has $406$ nonzero weights and is about $94\%$ sparse, obtaining an $R^2$-score of $99.9979$ on training data and $99.9977$ on test data. 
Figure~\ref{fig:adp} compares the input $\hat u_0$ of the ADP controller
(obtained by solving the convex approximation~\eqref{e-adp}, which uses the PCF $f$)
and the input $u^\star_0$ obtained by (approximately)
solving the nonlinear optimal control problem~\eqref{eq:optimal_control}
out to $t=T=150$,
for training and test data (generated similarly to training data).
We further compare the performance of the ADP controller to that of the
nonlinear controller,
for swinging up the pendulum from the initial state $z_0=(0, 0)$, and $\theta = m=1$.
The resulting trajectories are shown on the right of figure~\ref{fig:adp}.
Albeit the simplicity of the ADP controller, it behaves similarly to the nonlinear controller.
%$t$ from the current state $z_t$.

%The two controllers behave similarly, although the first only solves a convex optimization
%problem with one variable ($u_t$), while the second one requires solving a nonconvex
%nonlinear program with $H=150$ variables ($u_0,u_1,\ldots,u_{H-1}$).

%\subsection{Resource allocation}
%Notes: Fit concave utility where parameter could be context like time of day.
%
%\paragraph{Experimental setup.} 
%Tbd.
%
%\paragraph{Results.}
%Tbd.
%
%\subsection{\AB{Monotonicity example (utility function?)}}
%
%\subsection{\AB{Convex regression example}}

\clearpage
\section{Conclusions}\label{s-conclusions}
We have shown how to fit a PCF to data,
in a simple yet customizable way, with our open-source Python tool LPCF.
Our method allows (parametrized) convex optimization to be (in part)
data driven: while the modeling stage may heavily rely on
learning functions from data, the first-principles structure of convex
optimization is retained when solving a given problem instance.
Our experiments exhibit good modeling accuracy (also compared to
alternative convex function approximations), and the use of learned
PCFs in convex optimization problems.

%\subsection*{Acknowledgments}

\clearpage
{\small
	\bibliography{refs}
}

%\clearpage
%\appendix
%\section{Appendix}

\end{document}